\AddToHook{begindocument/before}{} 
\DeclareHookRule{begindocument}{hyperref}{before}{my-mdpi} 

\documentclass[preprints,article,accept,moreauthors,pdftex]{my-mdpi}

\firstpage{1} 
\pubvolume{12}
\issuenum{3}
\articlenumber{434}
\pubyear{2024}
\copyrightyear{2024}
\externaleditor{Academic Editor: Panayiotis Vafeas}
\datereceived{17 November 2023} 
\daterevised{22 January 2024} 
\dateaccepted{26 January 2024} 
\datepublished{29 January 2024} 
\hreflink{https://doi.org/10.3390/\linebreak math12030434} 
\doinum{10.3390/math12030434}
\pdfoutput=1 

\Title{Existence and Uniqueness of Weak Solutions to~Frictionless-Antiplane Contact Problems}

\TitleCitation{Existence and Uniqueness of Weak Solutions to Frictionless-Antiplane Contact Problems}

\Author{Besma Fadlia $^{1}$\orcidA{}, 
Mohamed Dalah $^{2}$\orcidB{} 
and Delfim F. M. Torres $^{3,4,}$*\orcidC{}}

\AuthorNames{Besma Fadlia, Mohamed Dalah and Delfim F. M. Torres}

\AuthorCitation{Fadlia, B.; Dalah, M.; Torres, D.F.M.}

\address{$^{1}$ \quad Laboratory of Differential Equations, Department of Mathematics,
University of Constantine 1,\newline 
Ain El Bey Road, {Constantine} P.O. Box 325, Algeria; 
besma.fadlia@student.umc.edu.dz\\  

$^{2}$ \quad Laboratory of Applied Mathematics and Modeling, 
Department of Mathematics, University of Constantine 1,
Ain El Bey Road, {Constantine} P.O. Box 325, Algeria; 
dalah.mohamed@umc.edu.dz \\
$^{3}$ \quad Center for Research and Development in Mathematics and Applications (CIDMA),\newline
Department of Mathematics, University of Aveiro, 3810-193 Aveiro, Portugal\\
$^{4}$ \quad Faculty of Business and Communications, 
INTI International University,\newline 
Persiaran Perdana BBN, Putra Nilai, 71800 Nilai, Negeri Sembilan, Malaysia}

\corres{Correspondence: delfim@ua.pt; Tel.: +351-234-370-668}

\abstract{We investigate a quasi-static-antiplane contact problem, 
examining a thermo-electro-visco-elastic material with 
a friction law dependent on the slip rate, assuming that 
the foundation is electrically conductive. The mechanical 
problem is represented by a system of partial differential equations, 
and establishing its solution involves several key steps. {Initially}, 
we obtain a variational formulation of the model, which comprises three systems: 
a hemivariational inequality, an elliptic equation, and a parabolic equation. 
{Subsequently}, we demonstrate the existence of a unique weak solution to the model. 
The proof relies on various arguments, including those related to evolutionary 
inequalities, techniques for decoupling unknowns, 
and certain results from differential equations.}
	


\begin{document}


\section{Introduction} 
\label{sect1}

Over the past few decades, phenomena related to the contact 
between deformable bodies have had a significant impact on society. 
The contact of brakes with wheels, a ball with the ground, 
and the application of a force to cutting objects, are just a few everyday 
instances, among many other examples. {Consequently}, research in this area is growing 
across various fields, particularly in engineering 
and mathematical literature~\cite{MR4083564,MR4130489}.

A comprehensive study of mechanical problems involves mathematical modeling. 
The modeling of {mechanical phenomena} is determined by a set of hypotheses 
that influence the system of partial differential equations. These hypotheses 
cover aspects such as the nature of the mechanical process 
(static, quasi-static, or dynamic), the behavior of the material 
(electromechanical, thermomechanical, elastic, etc.), and the boundary 
conditions on the contact surface (friction, adhesion, etc.)
\cite{MR4638641,MR4659364}.

In 1933, notable advancements in the mathematical and mechanical 
exploration of contact-mechanics problems took place. 
Signorini seems to have been the first to delve into this subject 
when he formulated the challenge of contact between a deformable body 
and a foundation. The solution to this problem eventually emerged in 1964, 
courtesy of Fichera, who employed emerging mathematical techniques. 
The initial findings concerning the existence and uniqueness of contact 
problems were established by Duvaut and Lions. Following this, a multitude 
of subsequent works have concentrated on resolving these variational 
problems \cite{a113,a12}.

Mathematics {plays} a pivotal role in the field of contact mechanics 
by providing diverse contact conditions. Contemporary research 
in this field also encompasses behavior laws that establish connections 
between mechanical and electrical effects, known as 
piezoelectricity\mbox{~\cite{MR4655630,MR4557105}}. This area is extensively 
explored in engineering structures, due to its 
temperature-dependent characteristics. Several articles have tackled 
the thermo-piezoelectric contact problem with friction
\cite{a16,a21,a22,a113,a23,a24}, 
thermo-elasto-visco-plasticity \cite{a15,aAZ,aEZ}, 
and thermo-visco-elasticity \cite{a13,a188}.
These studies are rooted in variations 
in constitutive laws and contact conditions 
\cite{25d9,25d10,au4,rk446,a255,25d11,q000}.

Here we propose a novel mathematical model designed to address 
the frictional-antiplane contact problem between a thermo-piezoelectric 
body and a conductive foundation. The innovation within our model lies 
in the introduction of a novel frictional contact tailored for such materials. 
This model incorporates a slip-rate-dependent friction law and introduces 
modifications to the electrical and thermal conditions used to describe the contact. 
We specifically focus on the case of antiplane-shear deformation.

The manuscript is organized as follows. In Section~\ref{sect2}, 
we recall some notations and definitions used in the sequel. 
In Section~\ref{sect3a}, we introduce our mathematical model 
for the quasi-static thermo-electro-visco-elastic-antiplane contact 
problem. Given some consistency assumptions about the data, 
the variational formulation of the model is then obtained 
in Section~\ref{sect3b}. Finally, in Section~\ref{sect4}, 
we prove {the} existence and uniqueness result 
{of} the formulated problem.


\section{Preliminaries} 
\label{sect2}

In this brief section, we review the fundamental notations 
and definitions that are employed throughout the work. 
For more  details, we refer the interested reader to \cite{q000}.

The space of second{-}order symmetric tensors on $\mathbb{R}^3$ 
is denoted by $\mathcal{S}^3$. The inner product and the  
norm on $\mathbb{R}^3$ are given as follows: 
$$
u^h \cdot v^h= u_i^h v_i^h,  
\quad \Vert  v^h \Vert=(v^h \cdot v^h)^\frac{1}{2} 
\quad \text{for all} \quad 
u^h=(u_i^h), v^h=(v_i^h)\in \mathbb{R}^3, 
\ 1\leq i,j\leq 3.
$$ 
The inner product on $\mathcal{S}^3$ is defined by
$$ 
\sigma^h \cdot \tau^h=\sigma_{i,j}^h\tau_{i,j}^h,  
\quad \text{for all} \quad 
\sigma^h=(\sigma_{i,j}^h), 
\quad \tau^h=\tau_{i,j}^h\in \mathcal{S}^3, 
\ 1\leq i,j\leq3,
$$
while its corresponding norm is given by
$$
\Vert \tau^h \Vert=\left(\tau^h \cdot \tau^h\right)^\frac{1}{2}.
$$
Next, we consider the following function spaces:
$$
V^h=\lbrace v^h \in H^{1}(\Omega^h), 
v^h=0 \ \text{on} \ \Gamma_{1}^h \rbrace,
$$
$$
W^h=\lbrace \psi^h \in H^{1}(\Omega^h), 
\psi^h=0 \ \text{on} \ \Gamma_{1}^h \rbrace ,
$$
$$
E^h=\lbrace \theta^h \in H^{1}\left( \Omega^h \right) ,\text{\ }\theta^h {=0}
\ \text{on} \ \Gamma^h_{1}\cup \Gamma_{2}^h \rbrace.
$$
It is an established fact that $V^h$, $W^h$, and $E^h$ 
are real Hilbert spaces with the inner products
$$
(u^h,v^h)_{V^h}=\int_{\Omega^h} \nabla u^h \cdot \nabla v^h  dx, 
\quad \forall u^h, v^h \in V^h,
$$
$$
(\varphi^h, \psi^h)_{W^h}
=\int_{\Omega^h} \nabla \varphi^h \cdot \nabla \psi^h  dx, 
\quad \forall \varphi^h, \psi^h \in W^h,
$$
$$
(\theta^h, \phi^h)_{E^h}
=\int_{\Omega^h} \nabla \theta^h \cdot \nabla \phi^h  dx, 
\quad \forall  \theta^h, \phi^h \in E^h.
$$
Furthermore, the associated norms 
\begin{align} 
\label{9}
\Vert v^h \Vert_{V^h} 
&=\Vert \nabla v^h \Vert_{L^{2}(\Omega^h)^{2}}, 
\quad \forall v^h \in V^h,\\  \label{10}
\Vert \psi^h \Vert_{W^h} 
&=\Vert \nabla \psi^h \Vert_{L^{2}(\Omega^h)^{2}}, 
\quad \forall \psi^h \in W^h,\\
\Vert \theta^h \Vert_{E^h} 
&=\Vert \nabla \theta^h \Vert_{L^{2}(\Omega^h)^{2}}, 
\quad \forall \theta^h \in E^h
\end{align}
on $V^h$, $W^h$, and $E^h$ are equivalent to the usual norm 
$\Vert \cdot \Vert_{H^{1}(\Omega^h)}$.
In the  light of Sobolev's trace theorem, we deduce that
\begin{align}
\label{11}
\exists \ c_{V^h} > 0 :  \Vert v^h \Vert_{L^{2}(\Gamma_{3}^h)} 
& \leq c_V^h \Vert v^h \Vert_{V^h}, 
\quad \forall v^h \in V^h, \\ \label{12}
\exists \  c_{W^h} > 0 : \Vert v^h \Vert_{L^{2}(\Gamma_{3}^h)} 
& \leq c_W^h \Vert \psi^h \Vert_{W^h}, 
\quad \forall \psi \in W^h,\\
\exists \  c_{E^h} > 0 : \Vert v^h \Vert_{L^{2}(\Gamma_{3}^h)} 
& \leq c_E^h \Vert \theta^h \Vert_{E^h}, 
\quad \forall \theta \in E^h.
\end{align}
In a real Banach space $(X^h, \Vert \cdot \Vert_{X^h})$, 
we denote the spaces of continuous and continuously differentiable 
functions on $[0, T]$ into $X^h$ by $C(0,T,X^h)$ and $C^{1}(0,T,X^h)$, 
respectively. These spaces are equipped with their respective norms,
$$
\Vert x^h \Vert_{C(0,T,X^h)}
= \max_{t\in [0, T]} \Vert x^h \Vert_{X^h}
$$
and 
$$
\Vert x^h \Vert_{C^{1}(0,T,X^h)}
= \max_{t\in [0, T]} \Vert x^h \Vert_{X^h}
+\max_{t\in [0, T]} \Vert \dot{x}^h \Vert_{X^h}.
$$
Below, we shall use the notations for the Lebesgue space $L^{2}(0, T, X^h)$ as well as
the Sobolev space $W^{1,2}(0, T, X^h)$. Recall that the  norm on  $L^{2}(0, T, X^h)$  
is given by
$$
\Vert u^h\Vert_{L^{2}(0, T, X^h)}^2=\int_{0}^T \Vert u^h(t)\Vert_{X^h}^{2}dt.
$$
Moreover,
$$
\Vert u^h\Vert_{W^{1,2}(0, T, X^h)}^2
=\int_{0}^T \Vert u^h(t)\Vert_{X^h}^{2}dt
+ \int_{0}^T \Vert \dot{u}^h(t)\Vert_{X^h}^{2}dt
$$
defines a norm on $W^{1, 2}(0, T, X^h)$.


\section{Formulation of the Problem}
\label{sect3a}

We begin by making explicit the physical setting under investigation.
Let us consider a thermo-electro-visco-elastic body $\mathcal{B}^{h}$ 
with friction. Assume that $\mathcal{B}^{h} \in \mathbb{R}^3$ is a cylinder 
with generators parallel to the $x_3^{h}$-axes having a regular region in 
its cross{-}section $\Omega^{h}$  in the $x_1^{h}, x_2^{h}$-plane, 
$Ox_1^{h}x_2^{h}x_3^{h}$ being a Cartesian-coordinate system. The cylinder 
is presumptively long enough for the end effects in the axial direction 
to be insignificant. Thus, $\mathcal{B}^{h}= \Omega^{h} \times ( -\infty, +\infty)$.
Let $\partial \Omega^{h} =\Gamma^{h}$. We suppose that $\Gamma^{h}$ is divided 
into three  disjoint measurable parts: $\Gamma_1^{h}$, $\Gamma_2^{h}$, and $\Gamma_3^{h}$. 
One has a partition of $\Gamma_1^{h} \cup \Gamma_2^{h}$ into 
two open parts, $\Gamma_a^{h}$ and $\Gamma_b^{h}$, such that the measures 
of $\Gamma_a^{h}$ and $\Gamma_b^{h}$ are positive. For $T > 0$, we denote 
by $[0, T]$ the time interval of interest. The $\mathcal{B}^{h}$ body is 
affected to time-dependent forces $\mathbf{f_0^{h}}$ and has a volume 
with free electric charges of density $q_0^{h}$. The body is clamped on 
$\Gamma_1^{h} \times (-\infty, +\infty)$ and, therefore, the displacement 
field disappears there. The surface tractions of density $\mathbf{f_2^{h}}$ 
act on $\Gamma_2^{h} \times (-\infty, +\infty)$. We also assume that the 
electrical potential vanishes on  $\Gamma_a^{h} \times (-\infty, +\infty)$ 
and that a surface electrical charge of density $q_2^{h}$ is prescribed on 
$\Gamma_b^{h} \times (-\infty, +\infty)$. The body is in contact  
with a rigid foundation $\Gamma_3^{h} \times (-\infty, +\infty)$ 
\cite{bou6,naaz}. We assume that
\begin{align}
\label{ad}
q_0^h
&=q_0^h(x_1^h,x_2^h,t):\Omega^h \times [0, T]
\longrightarrow \mathbb{R}, \\ \label{ae}  
q_2^h
&=q_2^h(x_1^h,x_2^h,t): \Gamma_b^h \times [0, T]
\longrightarrow \mathbb{R}, \\ \label{aa}
{\mathbf{f}_0^h}=(0,0,f_0^h),
& \  \text{with} \ f_0^h=f_0^h(x_1^h,x_2^h,t):\Omega^h \times [0,T]
\longrightarrow \mathbb{R}, \\ \label{ab}
{\mathbf{f}_2^h}=(0,0,f_2^h),
& \ \text{with} \  f_2^h=f_2^h(x_1^h,x_2^h,t):\Gamma_2^h \times [0,T]
\longrightarrow \mathbb{R}.
\end{align}
The electric charges (\ref{ad}) and (\ref{ae}) and the forces (\ref{aa}) and (\ref{ab}) 
are expected to give rise to electric charges and to deformations of the 
piezoelectric cylinder corresponding to an electric-potential field $\varphi^h$ 
and to a displacement ${\mathbf{u}^h}$,  which are independent of $x_3^h$ and have the form
\begin{align} 
\label{ar}
\varphi^h=
& \varphi^h(x_1^h,x_2^h,t):\Omega^h 
\times [0,T]\longrightarrow \mathbb{R}, \\ \label{ac}
{\mathbf{u}^h}=& u^h(0,0,u^h), \  \text{with} 
\ u^h=u^h(x_1^h,x_2^h,t):\Omega^h 
\times [0,T]\longrightarrow\mathbb{R}.
\end{align}
We denote by $\theta^h$ the temperature field, which is of the form
\begin{equation} 
\label{yh}
\theta^h= \theta^h(x_1^h,x_2^h,t): 
\Omega^h \times [0,T]\longrightarrow \mathbb{R},
\end{equation}
by $\varepsilon^h(u^h)= (\varepsilon_{i,j}^h(u^h))$ 
the infinitesimal strain tensor, by $\sigma^h=(\sigma_{ij}^h)$ 
the stress field, by $E^h(\varphi^h)=(E_i^h(\varphi^h))$ 
the electric field, and by $D^h=(D_i^h)$ the electric displacement 
field, where
\begin{align} 
\label{au}
\varepsilon_{i,j}^h(u^h)
&= \frac{1}{2}(u_{i,j}^h+u_{j,i}^h),\\ \label{ai}
E_i^h(\varphi^h)
&=-\varphi_{,i}^h.
\end{align}
The material is modeled by a thermo-electro-visco-elastic 
constitutive law of the form
\begin{align} 
\label{ao}
\sigma^h&=2\alpha^h \varepsilon^h(\dot{u}^h)+\zeta^h tr 
\varepsilon^h(\dot{u^h})I+2\mu^h\varepsilon^h(u^h)
+\lambda^h tr \varepsilon(u^h)I-\mathcal{E}^{h*} 
E^h(\varphi^h)-M_{e}^h\theta^h ,\\ \label{ap}
D^h&=\mathcal{E}^{h*}\varepsilon^h(u^h)
+\beta^h E^h(\varphi^h),
\end{align}
where $\zeta^h$ and $\alpha^h$ are viscosity coefficients, 
$\lambda^h$ and $\mu^h$ are the Lame coefficients, 
tr$\varepsilon^h(u^h)=\varepsilon_{ii}^h(u^h)$, $I$ is 
the unit tensor in $\mathbb{R}^3$, $\beta^h$ is the electric 
permittivity constant, $\mathcal{E}^h$ represents a 
third{-}order piezoelectric tensor, 
and $\mathcal{E}^{h*}$ is its transpose.
We assume that
\begin{equation} 
\label{aq}
\mathcal{E}^h\varepsilon^h
=
\begin{pmatrix}
e^h(\varepsilon_{13}^h+\varepsilon_{31}^h)\\
e^h(\varepsilon_{23}^h+\varepsilon_{32}^h)\\
e^h(\varepsilon_{33}^h)
\end{pmatrix},  
\quad \forall \varepsilon^h=(\varepsilon_{i,j}^h) \in \mathcal{S}^3,
\end{equation}
where $e^h$ is the piezoelectric coefficient. Additionally, 
$M_{e}^h= (M_{ij}^h)$ represents the thermal expansion tensor, 
which takes the form
\begin{equation}
M_{e}^h
=
\begin{pmatrix}
0 & 0 & M_{e_1}\\
0 & 0 & M_{e_2}\\
M_{e_1} & M_{e_2} & 0
\end{pmatrix}.
\end{equation}
Conversely, we assume that the coefficients $M_e^h$, $\alpha^h$, 
$\mu^h$, $\beta^h$, and $e^h$ depend on  $x_1^h$ and $x_2^h$. 
However, they are independent on $x_3^h$. Owing to $\mathcal{E}^h\varepsilon^h \cdot v^h
=\varepsilon^h \cdot \mathcal{E}^{h*}v^h$ for all $\varepsilon^h \in \mathcal{S}^3$
and $v^h \in \mathbb{R}^3$, it follows from $(\ref{aq})$ that
\begin{equation} 
\label{as}
\mathcal{E}^{h*}v^h
=
\begin{pmatrix}
0 & 0 & e^h v_{1}^h\\
0 & 0 & e^h v_{2}^h\\
e^h v_{1}^h & e^h v_{2}^h & e^h v_3^h
\end{pmatrix}, 
\quad \forall  \ v^h=(v_i^h) \in \mathbb{R}^3.
\end{equation}
In the antiplane context, having in mind (\ref{ar})--(\ref{yh}), 
and given the constitutive Equations (\ref{ao}) and (\ref{ap}) and equalities 
(\ref{aq}) and (\ref{as}), we obtain the stress field 
and the electric displacement field as 
\begin{equation} 
\label{aj}
\sigma^h
=
\begin{pmatrix}
0 & 0 & \sigma_{13}^h\\
0 & 0 & \sigma_{23}^h\\
\sigma_{31}^h & \sigma_{32}^h & 0
\end{pmatrix},
\ \mathbf{D}^h
=
\begin{pmatrix}
eu_{,1}^h-\beta^h\varphi_{,1}^h\\
eu_{,2}^h-\beta^h\varphi_{,2}^h\\
0
\end{pmatrix}, 
\end{equation} 
where
$$
\sigma_{13}^h=\sigma_{31}^h=\alpha^h \dot{u}_{,1}^h
+ \mu^h u_{,1}^h+e^h \varphi_{,1}^h-M_{e_{1}}\theta^h,
$$
$$
\sigma_{23}^h=\sigma_{32}^h=\alpha^{h} \dot{u}_{,2}^{h}
+\mu^{h}u_{,2}^{h}+e^{h}\varphi_{,2}^{h}-M_{e_{2}}\theta^h.
$$
We presume that the process is electrically static and mechanically 
quasi-static. Therefore, the equilibrium equations that govern it are
given by
\begin{equation}
\label{eq:ee}
\begin{aligned}
\mathrm{div}{ \sigma^h} + \mathbf{f}_0^h&=0, \\
\mathit{D}_{i,i}^h-q_0^h &=0 
\  \ \ \ \text{in} \ \mathcal{B}^h \times (0, T),
\end{aligned}
\end{equation}
where {$\mathrm{div}\sigma^h=(\sigma_{ij,j}^h)$} is the divergence 
of the tensor field {$\sigma^h$}. Accordingly, taking into 
account (\ref{ad}), (\ref{aa}), (\ref{ar})--(\ref{yh}), 
and  (\ref{aj}), the equilibrium Equation (\ref{eq:ee})
condenses into the following scalar equations:
\begin{align} 
\label{ak}
{\mathrm{div}}(\alpha^h\nabla \dot{u}^h
+\mu^h\nabla u^h)+{\mathrm{div}}(e^h\nabla \varphi^h)
-M_{e}^h {\mathrm{div}}(\theta^h)+f_{0}^h 
&= 0,    \   \ \text{in} \  \ \  \Omega^h \times(0, T), \\ \label{al}
{\mathrm{div}}(e^h\nabla u^h)-{\mathrm{div}}(\beta^h\nabla \varphi^h)
&=q_{0}^h,  \ \text{in} \  \Omega^h \times(0, T),
\end{align}
with 
\begin{equation} 
\label{lnj}
M_{e}^h
=
\begin{pmatrix}
M_{e_{1}}\\
M_{{e}_{2}}\\
0
\end{pmatrix}.
\end{equation}
In the sequel, we use the notation
$$
\mathrm{ div} {\tau^h}=\tau_{1,1}^h+\tau_{1,2}^h 
\ \ \text{for} \  
\tau^h=(\tau_1^h(x_1^h, x_2^h, t), \tau_2^h(x_1^h, x_2^h, t)),
$$
$$
\nabla v^h=(v_{,1}^h, v_{,2}^h), 
\ \partial_\nu^h v^h=v_{,1}^h\nu_{1}^h+v_{,2}^h\nu_2^h 
\ \  \text{for} \  v^h=v^h(x_1^h, x_2^h,t).
$$
Keeping in mind that the cylinder is clamped 
on $\Gamma_{1}^h \times(-\infty,+\infty)$, 
the electrical potential vanishes 
on $\Gamma_{a}^h \times(-\infty,+\infty)$,
and from (\ref{ad}) and (\ref{ac}) we find that
\begin{align} 
\label{am}
u^h=& 0   \ \ \  
\text{on}  \ \Gamma_{1}^h \times(0,T),\\ \label{aw}
\varphi^h=& 0 \ \ \ 
\text{on} \  \Gamma_{a}^h \times(0,T).
\end{align}
Note that $\nu^h $ is the unit normal on 
$\Gamma^h \times(-\infty,+\infty)$, where
\begin{equation} 
\label{ax}
\nu^h=(\nu_1^h, \nu_2^h,0), \text{with} \  
\nu_i^h=\nu_i^h(x_1^h,x_2^h): \Gamma^h 
\longrightarrow \mathbb{R}, \quad \text{for} \ i=1,2.
\end{equation}
We denote by $v_\nu^h$ and {$\mathbf{v}_\tau^h$} the normal 
and tangential components, respectively, 
of {$\mathbf{v}^h$} on the boundary---that is,
\begin{equation} 
\label{av}
v_\nu^h=\mathbf{v}^h \cdot \nu^h, 
\ \ \mathbf{v}_\tau^h={\mathbf{v}^h-v_\nu^h \nu^h}.
\end{equation}
Therefore, we denote by $\sigma_\nu^h$ and $\sigma_\tau^h$ 
the normal and the tangential components, respectively,  of 
$\sigma^h$ on the boundary, meaning that
\begin{equation}
\label{aù}
\sigma_\nu^h={(\sigma^h \nu^h) \cdot \nu^h},  
\ \ \  \sigma_\tau^h={\sigma^h\nu^h}
-\sigma_\nu^h {\nu^h}.
\end{equation}
From (\ref{aj}) and (\ref{ax}), we conclude that the Cauchy stress 
vector and the normal component of the electric displacement 
field are given by
\begin{equation} 
\label{an}
{\sigma^h\nu^h} =(0,0,\alpha^h \partial_{\nu^h} \dot{u}^h
+\mu^h \partial_{\nu^h} u^h+\partial_{\nu^h}\varphi^h-M_{e}^h \theta^h 
\cdot \nu^h), \ \ {D^h \cdot \nu^h}
=e^h \partial_{\nu^h} u^h-\beta^h \partial_{\nu^h} \varphi^h,
\end{equation}
respectively. Following this, we utilize the following notations:
$$
\partial_{\nu^h} u^h = u_{1}^h\nu_1^h+u_{2}^h\nu_2^h,
$$
$$
\partial_{\nu^h} \varphi^h 
= \varphi_{1}^h\nu_1^h+\varphi_{2}^h\nu_2^h. 
$$
Taking into account the traction boundary condition on 
$\Gamma_2^h \times(-\infty,+\infty)$ and the electric 
condition on $\Gamma_b^h \times(-\infty, +\infty)$, 
it follows from (\ref{ae}), (\ref{ab}), and (\ref{an}) that
\begin{align} 
\label{aaa}
\alpha^h \partial_{\nu^h} \dot{u}^h+\mu^h \partial_{\nu^h} u^h
+e^h\partial_{\nu^h}\varphi^h- M_{e}^h \theta^h \cdot \nu^h=f_2^h 
\ \ \text{on} \ \Gamma_2^h \times(0,T),\\ \label{aaz}
e^h \partial_{\nu^h} u^h-\beta^h \partial_{\nu^h} \varphi^h=q_2^h  
\ \ \text{on} \ \Gamma_b^h \times(0,T).
\end{align} 
We now describe the frictional contact condition on 
$\Gamma_3^h \times(-\infty,+\infty)$. First, we remark 
that from (\ref{ac}), (\ref{ax}), and (\ref{av}), we find 
that $u_\nu=0$, which indicates that the contact is bilateral. 
Thus, the contact is kept during the whole process. 
Now, using (\ref{ac}) and (\ref{ax})--(\ref{an}), we obtain
\begin{align}
\label{aae}
{u_\tau^h}&=(0,0,u^h), \\ \label{aar}
{\sigma_\tau^h}&=(0,0,\sigma_\tau^h),
\end{align}
where
\begin{equation} 
\label{aat}
\sigma_\tau^h = \alpha^h\partial_\nu^h \dot{u}^h
+\mu^h \partial_\nu^h u^h+e^h\partial_\nu^h
\varphi^h- M_{e}^h \theta^h \cdot \nu^h.
\end{equation}
We suppose that the friction is invariant {concerning} 
the $x_3^h$ axis, being modeled with a slip-rate-dependent 
friction law, where the strict inequality is satisfied 
in the stick zone and the equality in the slip zone---that is,
\begin{equation} 
\label{aay}
\begin{cases} 
\vert \sigma_\tau^h \vert\leq r(\vert \dot{u}_\tau^h\vert),  
\ \text{on}  \ \Gamma_3^h \times[0, T],\\
\sigma_\tau^h = -r(\vert \dot{u}^h\vert) 
\frac{ \dot{u}_\tau^h}{\vert \dot{u}_\tau^h\vert}, 
\  \   \text{on}  \ \Gamma_3^h \times[0, T].
\end{cases}
\end{equation}
Here, $r:\Gamma_3^h \times \mathbb{R} \longrightarrow \mathbb{R}^3$ 
is a given function, the friction bound, and $\dot{u}_\tau^h$ is the 
tangential velocity on the contact boundary: see \cite{au3,au8} for details. 
By (\ref{aae}) and (\ref{aat}), we ascertain that the conditions in (\ref{aay}) 
are given by
\begin{equation} 
\label{aau}
\begin{cases} 
\vert \alpha^h \partial_\nu^h \dot{u}^h+\mu^h \partial_\nu^h u^h
+e^h\partial_\nu^h\varphi^h- M_{e}^h \theta^h \cdot \nu^h \vert
\leq r(\vert \dot{u}^h\vert),\\
\alpha^h \partial_\nu^h \dot{u}^h+\mu^h \partial_\nu^h u^h
+e^h\partial_\nu^h\varphi^h- M_{e}^h \theta^h \cdot \nu^h
= -r(\vert \dot{u}^h\vert) \frac{ \dot{u}^h}{\vert \dot{u}^h \vert}, 
\   on  \ \Gamma_3^h \times[0, T].
\end{cases}
\end{equation}
Finally, we present the initial conditions:
\begin{equation}
u^h(0)=u_0^h,\quad \varphi^h(0)=\varphi_0^h, 
\quad \theta^h(0)=\theta_0^h 
\ \   \ \   \text{in} \ \Omega^h,  
\end{equation}   
where $u_0^h$, $\varphi_0^h$, and $\theta_0^h$ 
are given functions on $\Omega^h$.
Putting this all together, we obtain the problem
under investigation.

\begin{Problem}
\label{prb:01} 
Determine a displacement field 
$u^h:\Omega^h \times \left[ 0,T\right] \rightarrow 
\mathbb{R}$, the electric field
$\varphi^h:\Omega^h \times \left[ 0,T\right] \rightarrow 
\mathbb{R}$,  and a temperature field  
$\theta^h :\Omega^h \times \left[ 0,T\right] \rightarrow 
\mathbb{R}_+$, such that
\begin{align} 
\label{1}
{\mathrm{div}}(\alpha^h\nabla \dot{u}^h+\mu^h\nabla u^h)
+{\mathrm{div}}(e^h\nabla \varphi^h)
-M_e^h {\mathrm{div}}( \theta^h) +f_{0}^h 
&=0, \ \text{in} \ \Omega^h \times(0, T), \\ \label{2}
{\mathrm{div}}(e^h\nabla u^h)-{\mathrm{div}}(\beta^h\nabla \varphi^h)
&=q_{0}^h , \ \text{in} \ \Omega^h \times(0, T),\\ \label{3}
\dot{\theta}^h-{\mathrm{div}}(k\nabla \theta^h)=-M_e^h\nabla
&(\dot{u}^h)+p(t), \ \text{in} \ \Omega^h \times(0,T), \\ \label{313}
u^h&=0\ ,\ \text{on} \ \Gamma_1^h \times(0, T), \\ \label{4}
\alpha^h \partial_{\nu^h} \dot{u}^h+\mu^h \partial_{\nu^h}u^h
+e^h \partial_{\nu^h}\varphi^h-M_e^h \partial_{\nu^h}\theta^h 
&= f_{2}^h, \ \text{on} \ \Gamma_2^h \times(0, T),
\end{align}
\begin{equation} 
\label{5}
\begin{cases} 
\vert \alpha^h \partial_{\nu^h} \dot{u}^h+\mu^h \partial_{\nu^h} u^h
+\partial_\nu^h\varphi^h- M_{e}^h \theta^h.\nu^h \vert\leq r(\vert \dot{u}^h\vert),\\
\alpha^h \partial_{\nu^h} \dot{u}^h+\mu^h \partial_{\nu^h} u^h
+\partial_\nu^h\varphi^h - M_{e}^h \theta^h.\nu^h= -r(\vert \dot{u}^h\vert)  
\frac{ \dot{u}^h}{\vert \dot{u}^h\vert}, \  \text{on} \ \Gamma_3^h \times(0, T),
\end{cases}
\end{equation}
\begin{align}
\label{6}
\varphi^h&=0, \ \text{on} \  \Gamma_{a}^h \times(0, T), \\  \label{7}
e^h\partial_{\nu^h}u^h-\beta^h\partial_{\nu^h}\varphi^h
&=q_{2}^h, \ \text{on} \  \Gamma_b^h \times(0, T), \\  \label{8}
\theta_{0}^h&=0, \ \text{on} \ \Gamma_{1}^h
\cup \Gamma_{2}^h  \times(0, T), \\ \label{818}
-\mathcal{K}_{i,j}^h\frac{\partial \theta^h}{\partial x_j^h}v_i^h
&=\mathcal{K}_e^h(\theta^h-\theta_R^h), 
\ \text{on} \ \Gamma_{3}^h \times (0, T), \\ \label{888}
u^h(0)&=u_0^h, \ \quad \varphi^h(0)=\varphi_0^h, \  
\quad \theta^h(0)=\theta_0^h, \ \text{in} \ \Omega^h \times(0, T).
\end{align}
\end{Problem}

{It should be noted that the problem we consider here is different from that of \cite{25d10}.
Indeed, in \cite{25d10} they studied electro-visco-elastic material, meaning its
behavior was influenced solely by electrical and mechanical factors. The focus was purely on
the electrical impact without thermal interference, represented by equation (2.13)
of \cite{25d10}. By contrast, in our paper we add the factor of temperature, 
which plays a significant role in altering electromechanical properties 
through heating and cooling. For instance, the electro-rheological frictional 
force may vary depending on the temperature, thereby impacting the rheological 
properties of the material. Therefore, we can say that the material's behavior in 
our study is influenced by mechanical, electrical, and thermal
factors, represented by \mbox{Equation \eqref{ao}}, followed by the associated
conditions boundary in the Problem~\ref{prb:01}. 
On the other hand, the addition of the temperature factor also provides us with an
opportunity to explore and search for a different and simpler mathematical-solution
method (the decoupling-of-unknowns method). Additionally, there is a difference between 
the friction conditions. Indeed, in~\cite{25d10} the authors used Tresca's 
friction law (equation (2.25) there), which is characterized by setting 
a maximum limit for slip. It relies solely on the maximum force
without considering time rates. This was what the equation described there, observing that
the friction bound was constant and had one condition, making it easier to solve the
posed problem. By contrast, in our paper we use the slip-rate-dependent friction law \eqref{aay}:
here, the friction level depends on the rate of slip changes over time. This implies that
variations in the slip rate impact the friction level (the friction being higher when the
rate of change is greater and lower when the rate of change is smaller). This behavior is
described by Equation \eqref{aay}, where we represent friction bound by the variable $r$.}


\section{Variational Formulation}
\label{sect3b}
{To study} our Problem~\ref{prb:01}, we assume that the viscosity coefficient, 
the electric permittivity coefficient, the Lame coefficient, 
and the piezoelectric coefficient, satisfy 
\begin{align} 
\label{13}
&\alpha^h \in L^{\infty}(\Omega^h) \ \text{and there exists} 
\ \alpha^{h*} > 0, \text{ such that }\alpha^h(x^h) 
\geq \alpha^{h*}, \text{ a.e., } x^h \in \Omega^h,\\ \label{14}
&\beta^h \in L^{\infty}(\Omega^h) \ \text{and there exists} \ \beta^{h*} > 0, 
\text{ such that } \beta^h(x^h)\geq \beta^{h*}, 
\text{ a.e., } x^h \in \Omega^h, \\ \label{15}
&\mu^h \in L^{\infty} \ \text{and} \ \mu^h(x^h)> 0, 
\text{ a.e., }x^h \in \Omega^h,\\ \label{16}
&e^h \in  L^{\infty},
\end{align}
respectively. We also assume that the thermal tensors 
$M_e^h=\left( M_{ij}\right):\Omega^h
\rightarrow \mathcal{S}^3$ satisfy
\begin{equation}
\label{4.12}	
M_{ij}=M_{ji}\text{ }\in L^{\infty }\left( \Omega^h\right),
\text{\ }1\leqslant i,j\leqslant 3.  
\end{equation}
The thermal conductivity tensor 
$\widetilde{\mathcal{K}^h}=\left( \mathcal{K}_{ij}^h\right)
:\Omega^h \rightarrow \mathcal{S}^3$ satisfies
\begin{equation}
\left\{ 
\begin{array}{l}
\hbox{(i)}\text{ }\mathcal{K}_{ij}^h=\mathcal{K}_{ji}^h\text{ }
\in L^{\infty}\left( \Omega^h \right) ,\text{\ }1
\leqslant i,j\leqslant 3, \\ 
\hbox{(ii)}\text{ there exists }m_{\mathcal{K}^h}>0,\ \text{such that} \\ 
\mathcal{K}^h s\cdot s\geq m_{\mathcal{K}^h}\left\vert s\right\vert ^{2}, 
\text{ a.e., } x \in \Omega^h \text{, } \forall s\in \mathbb{R}^{3}.
\end{array}
\right.  \label{4.14}
\end{equation}
The boundary thermal data satisfy
\begin{equation}
\label{4.16}	
p^h\hspace{0.1cm}\in W^{1,2}(0,T;L^{2}(\Omega^h)),
\ \theta_{R}^h \in  W^{1,2}(0,T;L^{2}(\Gamma_3^h)), 
\  K_e^h \in L^{\infty }\left( \Omega^h, \mathbb{R_{+}}\right). 
\end{equation}
The forces, tractions, volume, and surface free charge 
densities have the following regularity: 
\begin{align} 
\label{17}
f_{0}^h \in W^{1,2}(0, T, L^{2}(\Omega^h)), 
\  f_{2}^h \in W^{1,2}(0, T,L^{2}(\Gamma_{2}^h)) , \\  \label{18}
q_0^h \in W^{1,2}(0, T, L^{2}(\Omega^h)), \\ \label{19}
q_2^h \in W^{1,2}(0, T, L^{2}(\Gamma_b^h)), 
\  q_2^h =0,  \text{ a.e.,} \ x^h \in \Gamma_b^h,
\end{align}
respectively. The friction bound satisfies
\begin{equation} 
\label{20}
\begin{cases} 
a)\  r:\Gamma_3^h \times \mathbb{R} \longrightarrow \mathbb{R^+},\\
b)\  \exists \  L_r > 0,  \ \text{ such that} 
\  \vert r(x^h, s_1)-r(x^h, s_2) \vert \leq L_r\vert s_1 -s_2 \vert,\\
c)\  x^h \longrightarrow r(x^h,s) \ \text{is Lebesgue-measurable on} 
\  \Gamma_3^h, \ \forall \ s \in \mathbb{R},\\
d)\ \text{the mapping } x^h \longrightarrow r(x^h, 0) 
\ \text{belongs to} \ L^{2}(\Gamma_3^h).
\end{cases}
\end{equation}
We define the functional 
$j:v^h \times v^h \longrightarrow \mathbb{R} $ by
\begin{equation} 
\label{22}
j(\dot{u}^h,v^h)= \int_{\Gamma_3^h} r(\vert \dot{u}^h \vert)
\vert v^h \vert da, \ \forall v^h \in V^h.
\end{equation}
Let $\eta_1$, $\eta_2$, $v_1$, $v_2 \in V$. 
By using (\ref{20}) and (\ref{22}), we find that
$$
j(\eta_{1}, v_{2}^h)-j(\eta_{1}, v_{1}^h)+j(\eta_{2}, v_{1}^h)
-j(\eta_{2}, v_{2}^h) \vert \leq L_r  \Vert \eta_{1} 
- \eta_{2} \Vert_{L^2(\Gamma_3^h)}  \Vert  v_{1} 
- v_{2} \Vert_{L^2(\Gamma_3^h)}, 
$$
{and applying the norm mentioned in (\ref{11}), one obtains
$$ L_r  \Vert \eta_{1} 
- \eta_{2} \Vert_{L^2(\Gamma_3^h)}  \Vert  v_{1} 
- v_{2} \Vert_{L^2(\Gamma_3^h)} \leq c_V^2 L_r  \Vert \eta_{1} 
- \eta_{2} \Vert_{V^h}  \Vert  v_{1}^h - v_{2}^h \Vert_{V^h}.$$
Hence, we conclude that }
\begin{align} 
\label{23}
j(\eta_{1}, v_{2}^h)-j(\eta_{1}, v_{1}^h)+j(\eta_{2}, v_{1}^h)
-j(\eta_{2}, v_{2}^h)  \leq c_V^2 L_r  \Vert \eta_{1} 
- \eta_{2} \Vert_{V^h}  \Vert  v_{1}^h - v_{2}^h \Vert_{V^h}.
\end{align}
The initial data verify
\begin{equation} 
\label{21}
u_0^h \in V^h , \  \theta_0^h \in L^2(\Omega^h).
\end{equation}
We use functions $f^h:[0,T] \longrightarrow V^h$ 
and $ q^h:[0, T] \longrightarrow W^h $ as
\begin{align} 
\label{24}
(f^h,v^h)_{V^h} = \int_{\Omega^h} f_{0}^h v^h dx^h 
+ \int_{\Gamma_{2}^h}f_{2}^h v^h da,  \    \forall v^h \in V^h,\\ \label{25}
(q^h,\psi^h)_{V^h} = \int_{\Gamma_b^h }q_{2}^h \psi^h dx^h 
+ \int_{\Omega^h} q_{0}^h \psi^h  da,  \  \forall \psi^h \in W^h.
\end{align}
The definitions of $f^h$ and $q^h$ are based on Riesz's representation theorem. 
Therefore, by using assumptions (\ref{24}) and (\ref{25}), we ascertain that 
the  above integrals are well defined and that
\begin{align} 
\label{26}
f^h \in  W^{1,2}(0, T, V^h), \\ \label{27}
q^h \in W^{1,2}(0, T, W^h).
\end{align}
Then, the function $P:[0, T]\longrightarrow (E^{h})^{\prime}$ 
and the operators $\widetilde{\mathcal{K}^h}:E \longrightarrow (E^{h})^{\prime},
\widetilde{\mathcal{M}^h}:V^h \longrightarrow (E^{h})^{\prime}$ are defined by
$$
\langle P(t),\mu^h\rangle_{(E^{h})^{\prime}  
\times E^h} =\displaystyle \int_{\Omega^h}p^h \mu^h dx^h 
+ \displaystyle \int_{\Gamma_{3}^h}\mathcal{K}_e^h \theta_{R}^h \mu^h ds,
$$
$$
\langle \widetilde{\mathcal{K}^h}\tau,\mu^h\rangle_{(E^{h})^{\prime}  
\times E^h} =\displaystyle \int_{\Gamma_{3}^h}\mathcal{K}_{e}^h\tau\mu^h ds 
+ \sum^{d}_{i,j=1} \displaystyle \int_{\Omega^h}
\mathcal{K}_{i,j}^h\frac{\partial \mu^h}{\partial x_{j}^h} 
\frac{\partial \mu^h}{\partial x_{i}^h}dx,
$$
$$
\langle \widetilde{\mathcal{M}^h}v^h,\mu^h\rangle_{(E^{h})^{\prime}  
\times E^h} = -\int_{\Omega^h}(M_{e}^h \nabla v^h )\mu^h dx 
+\int_{\Gamma_{3}^h}h_{\tau}(\vert v_{\tau}^h\vert)\mu^h ds 
$$
for all $v^h \in V^h, \ \tau \in E^h, \ \mu^h \in E^h$. 
Moreover, we define the bilinear forms 
$a_{\alpha^h}: V^h \times V^h \longrightarrow \mathbb{R}$, 
$a_{\mu^h}: V^h\times V^h \longrightarrow \mathbb{R}$, 
$a_{e^h}:V^h \times W^h \longrightarrow \mathbb{R}$, 
$a_{e^h}: W^h \times V^h \longrightarrow \mathbb{R}$, 
$a_{\beta^h}: W^h \times W^h \longrightarrow \mathbb{R}$ 
and $a_{M^h}: E^h \times V^h \longrightarrow \mathbb{R}$ by
\begin{align}
\label{28}
a_{\alpha^h}(u^h,v^h)
&= \displaystyle \int_{\Omega^h} \alpha^h \nabla u^h 
\cdot \nabla v^h dx,\\ \label{30}
a_{\mu^h}(u^h,v^h)
&= \displaystyle \int_{\Omega^h} \mu^h \nabla u^h
\cdot \nabla v^h dx,\\ \label{301}
a_{M^h}(u^h,v^h)
&=- \displaystyle \int_{\Omega^h} M_{e}^h  u^h
\cdot \nabla v^hdx ,\\ \label{29}
a_{e^h}(u^h, \varphi^h)
&= \displaystyle \int_{\Omega^h} e^h \nabla u^h
\cdot \nabla v^h dx = a_{e^h}(v^h,u^h), \\ \label{31}
a_{\beta^h}(\varphi^h,\psi^h)
&= \displaystyle \int_{\Omega^h} \beta^h \nabla \varphi^h 
\cdot \nabla \psi^h dx  
\end{align} 
for all $u^h, v^h \in V^h$, $\varphi^h, \psi^h \in W^h$. We note that 
by assumptions \eqref{13}--\eqref{16}, the  above integrals 
are well defined. {Using the definition of norms} 
\eqref{9}--\eqref{12}, we see that the forms $a_{\alpha^h}$, $a_{\mu^h}$, 
$a_{e^h} $, and $a_{\beta^h}$, are continuous and that the forms 
$a_{\alpha^h}$, $a_{\mu^h}$, and $a_{\beta^h}$ are symmetric. 
Furthermore, the form $a_{\alpha^h}$ is $V$-elliptic, i.e., 
\begin{align} 
\label{32}
a_{\alpha^h}(u^h,v^h) 
&\leq  \Vert \alpha^h \Vert_{L^{\infty}(\Omega^h)} 
\Vert u^h \Vert_{V^h}\Vert v^h \Vert_{V^h},  
\ \forall u^h, v^h \in V^h,\\ \label{33}
a_\alpha^h(v^h,v^h) 
&\geq \alpha^{h*} \Vert v^h \Vert_{V^h}^{2} 
\ \     \forall v^h \in V^h.
\end{align}

Now, we can state the variational formulation of our Problem~\ref{prb:01}.

\begin{Problem}[Variational Problem] 
Determine a displacement field $u^h:[0, T] \longrightarrow V^h$, 
an electric-potential field $\varphi^h:[0, T]\longrightarrow W^h$, 
and a temperature field $ \theta^h :[0, T] \longrightarrow E^h$, such that
\vspace{-14pt}

\begin{adjustwidth}{-\extralength}{0cm}
\begin{align} 
a_{\alpha^h}(\dot{u}^h(t),v^h-\dot{u}^h(t))+a_{\mu^h}(u^h(t),v^h-\dot{u}^h(t))
&+a_{e^h}(\varphi(t), v-\dot{u}^h(t)) 
+a_{M^h}(\theta^h(t), v^h-\dot{u}^h(t))\notag\\ \label{38}
+j(\dot{u}^h(t),v^h)-j(\dot{u}^h(t), \dot{u}^h(t)) 
&\geq (f^h,v^h-\dot{u}^h)_{V^h}, 
\ \forall v^h \in V^h,  \ t \in [0, T] ,\\ \label{39}
a_{\beta^h}(\varphi^h,\psi^h)-a_{e}(u^h,\psi^h)
&=(q^h,\psi)_W^h, \  \forall \psi^h \in W^h ,\\ \label{82}
\dot{\theta}^h(t)+\widetilde{\mathcal{K}^h}\theta^h(t)
&=\widetilde{\mathcal{M}^h}\dot{u}^h(t)+P(t), \text{ in } E^{h}, \\ \label{40}
u^h(0)&=u_0^h, \quad \varphi^h (0)=\varphi_0^h, \quad \theta^h(0)=\theta_0^h.
\end{align}
\end{adjustwidth}
\end{Problem}


\section{Existence and Uniqueness Result} 
\label{sect4}

In this section, we use the variational formulation 
to prove the existence and uniqueness of a weak solution to Problem~\ref{prb:01}.
For that, we make use of some auxiliary problems and lemmas. 

\begin{Problem}[Auxiliary Problem]
\label{tp1} 
{Determine} a displacement field 
$u^h: [0, T] \longrightarrow V^h$, such that
\vspace{-14pt}

\begin{adjustwidth}{-\extralength}{0cm}
\begin{align} 
\label{43}
a^h(u^h(t),v^h-\dot{u}^h(t))+ b^h(\dot{u}^h(t), v^h-\dot{u}^h(t))
&+j^h(\dot{u}^h(t),v^h)-j^h(\dot{u}^h(t), \dot{u}^h(t))\\ \notag
&\geq (F(t),v^h-\dot{u}^h(t))_{V^h}, \  \forall v^h \in V^h,  
\quad t \in [0, T], \\ \label{44}
u^h(0)&=u_0^h.
\end{align}
\end{adjustwidth}
\end{Problem}

In order to study Auxiliary Problem~\ref{tp1}, we assume that
\begin{equation} 
\label{45}
\begin{cases}
a^h:V^h \times V^h \longrightarrow \mathbb{R} 
\text{ is a bilinear form and there exists} \ M^h >0, \ \text{such that}  \\ 
\vert a^h(u^h,v^h) \vert \leq M^h \Vert u^h \Vert_{V^h} \Vert 
v^h \Vert_{V^h},  \   \forall u^h, v^h \in V^h;
\end{cases}
\end{equation}

\begin{adjustwidth}{-\extralength}{0cm}
\begin{equation} 
\label{46}
\begin{cases}
b^h: V^h \times V^h \longrightarrow 
\text{is a bilinear symmetric form verifying:} \\
(a)\ \text{There exists} \ M^{h'} >0, \  \text{such that }
\vert b^h(u^h,v^h) \vert \leq  \ M^{h'} \Vert u^h \Vert_{V^h} \Vert 
v^h \Vert_{V^h},  \ \forall u^h, v^h \in X^h, \\ 
(b)\ \text{there exists} \ m^{h'}> 0, \text{ such that }  
b^h(v^h, v^h)\geq m^{h'} \Vert v^h \Vert_{V^h}^2, 
\    \forall v^h \in V^h;
\end{cases}
\end{equation}
\end{adjustwidth}
\begin{equation} 
\label{47}
\begin{cases} 
j^h:V^h \times V^h \longrightarrow \mathbb{R} \  \text{satisfies:}\\
(a) \ \text{For all} \  \eta \in V^h,  \ j^h(\eta ,\cdot) 
\text{ is convex and 	l.s.c. on} \  V^h, \\
(b) \  \text{there exists} \ \alpha^{h'} \geq 0, 
\  \text{such that for all} \  \eta_1, \eta_2, v_1^h, v_2^h \in V^h 
\text{ we have}   \\  \vert j^{h}(\eta_{1}, v_{2}^h)-j^h(\eta_{1}, v_{1}^h)
+j^h(\eta_{2}, v_{1}^h)-j^h(\eta_{2}^h, v_{2}^h) \vert \leq \alpha^{h'}  \Vert 
\eta_{1} - \eta_{2} \Vert_{V^h}  \Vert  v_{1}^h - v_{2}^h \Vert_{V^h};
\end{cases}
\end{equation}
\begin{equation} 
\label{48}
u_0^h \in V^h;
\end{equation}
\begin{equation}
\label{49}
F \in W^{1,2}(0, T, V^h).
\end{equation}

From this point onward, we denote a generic constant as $c^h > 0$.

\begin{Lemma} 
\label{lp1}
Assume that \eqref{45}--(\ref{49}) holds. If $m^{h'} > \alpha^{h'}$, then
there exists a unique solution $u^h \in W^{1,2}(0, T, V^h)$ 
to Auxiliary Problem~\ref{tp1}.
\end{Lemma}

\begin{proof}
By using (\ref{43}), for any $t_{1}$ and $t_{2} \in [0, T]$ we find that
\vspace{-18pt}

\begin{adjustwidth}{-\extralength}{0cm}
\begin{align}  
\notag
a^h(u^h(t_{1}), v^h-\dot{u}^h(t_{1}))
&+ b^h(\dot{u}^h(t_{1}), v^h-\dot{u}^h(t_{1}))
+j^h(\dot{u}^h(t_{1}),v^h)-j^h(\dot{u}^h(t_{1}), \dot{u}^h(t_{1})) \\ \label{A12}
&\geq (F(t_{1}),v^h-\dot{u}^h(t_{1}))_{V^h}, 
\ \forall v^h \in V^h,\\ \notag \label{A13}
a^h(u^h(t_{2}), v^h-\dot{u}^h(t_{2}))
&+ b^h(\dot{u}^h(t_{2}), v^h-\dot{u}^h(t_{2}))
+j^h(\dot{u}(t_{2}),v)-j^h(\dot{u}^h(t_{2}), \dot{u}^h(t_{2}))\\
& \geq (F(t_{2}),v^h-\dot{u}^h(t_{2}))_{V^h},\ \forall v^h \in V^h.
\end{align}
\end{adjustwidth}
We take $v^h=\dot{u}^h(t_2)$ in (\ref{A12}) and 
$v^h=\dot{u}^h(t_1)$ in  (\ref{A13}) and, adding the two inequalities, we obtain
\vspace{-24pt}

\begin{adjustwidth}{-\extralength}{0cm}
\begin{align} 
\notag
&b^h(\dot{u}^h(t_1)-\dot{u}^h(t_2),\dot{u}^h(t_1)
-\dot{u}^h(t_2))+j(\dot{u}^h(t_1), \dot{u}^h(t_2))
-j^h(\dot{u}^h(t_1), \dot{u}^h(t_2))
+j^h(\dot{u}^h(t_2), \dot{u}^h(t_1))\\ 
&-j^h(\dot{u}^h(t_2), \dot{u}^h(t_2)) \leq a^h(u^h(t_1)- u^h(t_2), 
\dot{u}^h(t_2)- \dot{u}^h(t_1))+(F(t_1)-F(t_2),\dot{u}^h(t_2)-\dot{u}^h(t_1)).
\end{align}
\end{adjustwidth}
We now use assumptions (\ref{45})--(\ref{47}) to obtain
\begin{equation}
\Vert \dot{u}^h(t_1)-\dot{u}^h(t_2) \Vert_{V^h} 
\leq c^h( \Vert u^h(t_1)-u^h(t_2)\Vert_{V^h} 
+ \Vert F(t_1)-F(t_2) \Vert_{V^h}),
\end{equation}
where $c^h=\max \left\{ \displaystyle\frac{M^h}{m^{h'}
-\alpha^{h'}}, \frac{1}{m^{h'}-\alpha^{h'}}\right\}$. 
This inequality, combined with the regularity  
$u^h \in C^{1}(0, t, V^h)$, shows that 
$\dot{u}^h:[0, T]\longrightarrow V^h$ 
is an absolutely continuous function. 
Furthermore, it should be noted that
$$
\Vert \ddot{u}^h(t) \Vert_{V^h} \leq c^h( \Vert \dot{u}^h(t) \Vert_{V^h} 
+\Vert F(t)\Vert_{V^h}) \text{ a.e.,} \ t \in [0, T].
$$
We conclude that $u^h \in  W^{2,2}(0,T,V^h)$.
\end{proof}

Using the solution $u^h$ asserted by Lemma~\ref{lp1},
we proceed by considering a second auxiliary problem.

\begin{Problem}[Auxiliary Problem] 
\label{tp2}
Determine a function $\theta_\eta^h :[0 , T] \longrightarrow E^h$, such that
\begin{equation}
\left\{ 
\begin{array}{c}
\dot{\theta}_{\eta}^h(t)+\widetilde{\mathcal{K}^h}
\theta_{\eta}^h(t)=P(t)+\widetilde{\mathcal{M}^h}\dot{u}_\eta^h(t),\\
\theta_{\eta}^h(0)=\theta_0^h , \ \ t \in[0, T].
\end{array}
\right. \label{87}
\end{equation}
\end{Problem}

\begin{Lemma} 
\label{lp2}
For all $\eta \in C(0, T, V^h)$, there exists a unique solution 
to Auxiliary Problem~\ref{tp2} with
\begin{equation}
\label{88}	
\theta_\eta^h \in L^{2}\left( 0,T,E^h\right) \cap C\left( \left[ 0,T\right],
L^{2}\left( \Omega^h \right) \right) \cap W^{1,2}(0,T,(E^{h})^{\prime}).
\end{equation} 
Moreover, there exists $c^h > 0$, such that, 
for all $ \eta_1, \eta_2 \in C(0, T, V^h)$, we have
\begin{equation} 
\label{89}
\Vert \theta_{\eta_1}^h(t)-\theta_{\eta_2}^h(t)\Vert_{E^h} 
\leq c^h \int_{0}^{t} \Vert \dot{u}_{\eta_{1}}^h(s)
-\dot{u}_{\eta_{2}}^h(s)\Vert_{V^h}^2  \ ds, 
\quad \forall t \in [0, T].
\end{equation} 
\end{Lemma}

\begin{proof}
The existence and uniqueness result to (\ref{87})
follows by considering the Gelfand evolution
(see, e.g., \cite{MR1033497}):
$$
E^h \subset F^h \equiv (F^{h})^{\prime} \subset (E^{h})^{\prime},
$$
and verifying that the operator $\widetilde{\mathcal{K}^h}: E^h 
\longrightarrow (E^{h})^{\prime}$ is strongly monotonic 
and linearly continuous. As inferred from the expression 
of the operator $\widetilde{\mathcal{M}^h}$, we have
$$
\dot{u}_\eta^h(t) \in W^{1,2}(0,T;V^h) 
\ \Longrightarrow \widetilde{\mathcal{M}^h} 
\dot{u}_\eta^h(t) \in  W^{1,2}(0,T,F^h),
$$
and
$$
P(t) \in W^{1,2}(0,T,E^h) \ \Longrightarrow 
\widetilde{\mathcal{M}^h} \dot{u}_\eta^h(t)+P(t) 
\in W^{1,2}(0,T, (E^{h})^{\prime}).
$$
For $\eta_1, \eta_2  \in C(0, T, V^h)$, we have 
\begin{align}
& \left( \dot\theta_{\eta_1}^h(t)- \dot{\theta_{\eta_2}}^h(t),  
\theta_{\eta_1}^h(t) - \theta_{\eta_2}^h(t)\right) _{(E^{h})^{\prime}
\times E^h} +\left( \widetilde{\mathcal{K}^h}\theta_{\eta_1}^h(t)
- \widetilde{\mathcal{K}^h}\theta_{\eta_2}^h(t),
\theta_{\eta_1}^h(t)- \theta_{\eta_2}^h(t) \right)_{E^h} \\ \notag
&=\left( \widetilde{\mathcal{M}^h}\dot{u}_{\eta_1}^h(t)
-\widetilde{\mathcal{M}^h}\dot{u}_{\eta_2}^h(t),
\theta_{\eta_{1}}^h(t)- \theta_{\eta_2}^h(t) 
\right)_{L^{2}(\Omega^h )}, \quad t \in [0, T],
\end{align}
and, by integrating the last property over $(0, t)$ 
and utilizing the Lipschitz continuity of 
$\widetilde{\mathcal{M}^h}: V^h \longrightarrow (E^{h})^{\prime}$, 
as well as the strong monotonicity of $\widetilde{\mathcal{K}^h}$, 
we deduce that (\ref{89}) holds for $t \in [0, T]$.
\end{proof}

We are now ready to prove our main result: 
the existence of a unique solution to our Problem~\ref{prb:01}. 

\begin{Theorem} 
\label{t1}
Assume that (\ref{13})--(\ref{27}) are satisfied. Then, there exists $Z_0$, 
which depends on $\Omega^h$, $\Gamma_1^h$, $\Gamma_2^h$, and $\Gamma_3^h$, 
such that if $L_r <Z_0$ then there exists a unique solution 
$(u^h, \varphi^h, \theta^h)$ to (VP), satisfying
\begin{align}
\label{41}
u^h \in 
&  W^{2,2}(0, T, V^h) , \\ \label{132}
\theta^h \in  
& L^{2}\left( 0,T,E^h\right) \cap C\left( \left[ 0,T\right],
L^{2}\left( \Omega^h \right) \right) 
\cap W^{1,2}(0,T,(E^{h})^{\prime}),\\ \label{152}
\varphi^h \in  
& W^{1,2}(0,T,W^h).
\end{align}
\end{Theorem}

\begin{proof}
From (\ref{39}) we ascertain that
\begin{equation} 
\label{443}
(\beta^h \varphi^h, \psi^h)_{W^h}
-(e^h u^h, {\psi^h})_{W^h}=(q^h,\psi^h)_{W^h},
\end{equation}
while the use of (\ref{443}) gives 
$$
\beta^h \varphi^h(t)=e^hu^h(t)+q^h.
$$
Hence, we deduce that
\begin{equation} 
\label{444}
\varphi^h(t)=\dfrac{e^h}{\beta^h}u^h(t)
+\dfrac{q^h}{\beta^h}.
\end{equation}
From (\ref{82}) and (\ref{40}), we obtain the solution to
\begin{equation} 
\label{4444}
\theta^h(t)= \displaystyle \int_{0}^{t}\widetilde{\mathcal{M}^h}
e^{-\int_{s}^{t}\widetilde{\mathcal{K}^h} dx} \dot{u}^h(s)ds
+\displaystyle \int_{0}^{t} e^{-\int_{s}^{t}\widetilde{\mathcal{K}^h} dx}
P(s)ds+\theta_0^h e^{-\int_{0}^{t}\widetilde{\mathcal{K}^h} ds}.
\end{equation}
Taking (\ref{444}) and (\ref{4444}) and substituting them into (\ref{38}), we obtain
\begin{small}
\begin{align*}
a_{\theta^h}(\dot{u}^h(t), v^h-\dot{u}^h(t))+ a_{\mu^h}(u^h(t),v^h-\dot{u}^h(t))
+&a_{e^h}\left(\dfrac{e^h}{\beta^h}u^h(t)+\dfrac{q^h}{\beta^h}, v^h-\dot{u^h}(t)\right)\\
+a_{M^h}\Biggl(\int_{0}^{t}\widetilde{\mathcal{M}^h}e^{-\int_{s}^{t}
\widetilde{\mathcal{K}^h} dx} \dot{u}^h(s)ds
&+\int_{0}^{t} e^{-\int_{s}^{t}\widetilde{\mathcal{K}^h} dx}P(s)ds
+\theta_0^h e^{-\int_{0}^{t}\widetilde{\mathcal{K}^h} ds}, v^h-\dot{u}^h(t)\Biggr)\\
+j^h(\dot{u}^h(t),v^h)-j^h(\dot{u}^h(t), \dot{u}^h(t))
& \geq (f^h(t),v^h-\dot{u}^h)_{V^h}, 
\quad \forall v^h \in V^h, \ t \in [0, T], \\
u^h(0)&=u_0^h.
\end{align*}
\end{small}
Next, we define the bilinear forms $a^h: V^h \times V^h \longrightarrow \mathbb{R}$ 
and $b^h: V^h \times V^h \longrightarrow \mathbb{R}$ as follows:
\begin{align} 
\label{50}
a^h(u^h(t),v^h-\dot{u}^h(t))
&=a_{\mu^h}(u^h(t),v^h-\dot{u}^h(t))
+a_{e^h}\left(\frac{e^h}{\beta^h}u^h(t), v^h-\dot{u}^h(t)\right), \\ \notag
&+a_{M^h}\left(\int_{0}^{t}\widetilde{\mathcal{M}^h}e^{-\int_{s}^{t}
\widetilde{\mathcal{K}^h} dx} \dot{u}^h(s)ds, v^h-\dot{u}^h(t)\right),\\ \label{51}
b^h(\dot{u}^h(t), v^h-\dot{u}^h(t))
&=a_{\theta^h}(\dot{u}^h(t), v^h-\dot{u}^h(t)).
\end{align} 
We also consider the function $F:[0, T] \longrightarrow V^h$, defined by
\begin{align}\label{52}
(F(t),v^h-\dot{u}^h)_{V^h}
=&(f^h(t),v^h-\dot{u}^h)_{V^h}-a_{e^h}
\left(\dfrac{q^h}{\beta^h}, v^h-\dot{u}^h(t)\right)\\ \notag
&-a_{M^h}\left(\displaystyle \int_{0}^{t} e^{-\int_{s}^{t}
\widetilde{\mathcal{K}^h} dx}P(s)ds+\theta_0^h 
e^{-\int_{0}^{t}\widetilde{\mathcal{K}^h} ds}, v^h-\dot{u}^h(t)\right).
\end{align}
From (\ref{50})--(\ref{52}), we obtain the following results:
\begin{itemize}
\item the bilinear form $a^h$ and the initial data 
$u_0^h$ satisfy conditions (\ref{45}) and (\ref{48});

\item $f^h \in W^{1,2}(0,T, V^h)$, 
$q^h \in W^{1,2}(0,T, W^h)$, and $P \in W^{1,2}(0,T, E^h)$, 
combined with the definition of $F(\cdot)$ in (\ref{52}), 
satisfy (\ref{49});

\item for all $\eta \in V^h$, the functional 
$j^h(\eta, \cdot): V^h \longrightarrow \mathbb{R}$ 
is a continuous seminorm on $V^h$, thus satisfying the 
condition (\ref{47})-(a) (recall that $j^h$ also fulfills inequality 
(\ref{23}), demonstrating that condition (\ref{47})-(b) 
holds with $\alpha^{h'} = c_V L_r$);

\item from (\ref{32}), it is evident that the bilinear form $b^h$ 
satisfies condition (\ref{47}) with $m^{h'} = \alpha^{h*}$.
\end{itemize}
Now, we choose $Z_0 = \displaystyle\frac{\alpha^{h*}}{c_V^2}$, 
which depends on $\Omega^h$, $\Gamma_1^h$, $\Gamma_2^h$, $\Gamma_3^h$, 
and $\alpha^{h'}$. Then, if $L_r < Z_0$, we have $m^{h'} > \alpha^{h'}$. 
Therefore, (\ref{41}) and (\ref{132})  are a direct consequence 
of Lemmas~\ref{lp1} and \ref{lp2}. For the proof of (\ref{152}), 
let $\varphi^h : [0,T] \longrightarrow W^h$ 
be the electrical-potential field defined by (\ref{444}). We have
$$
\begin{cases} 
u^h \in W^{2, 2}(0, T, V^h)\\
q^h \in W^{1,2}(0, T, W^h)
\end{cases}
\Rightarrow \varphi^h \in W^{1,2}(0, T, W^h),
$$
which completes the proof.
\end{proof}


\section{Conclusions} 
\label{sect5}

{In this paper, we theoretically investigated a quasi-static-antiplane 
contact problem with a slip-rate-dependent friction law involving a 
thermo-electro-visco-elastic body.
The slip-rate-dependent friction law was interesting to investigate because
the friction level depends on the rate at which the slip changes over time. This 
implies that variations in the slip rate impact the friction level, 
with the friction force increasing as the slip-rate speed increases, 
and decreasing as the slip-rate speed decreases. This is more reasonable 
because it considers both the quasi-static equilibrium
conditions and the localized dynamic interactions that may occur within the
contact interface.} 

{We used Green's formula to derive 
the variational formulation of our problem. Following that, 
we demonstrated the existence and uniqueness of the weak solution, 
using various techniques, including the time-dependent variational 
equation, the variational-evolution equation, the differential equation, 
and the decoupling of unknowns. The advantages of the developed model
include both the quasi-static equilibrium conditions and the 
localized dynamic interactions that may occur within the contact 
interface. This approach provides a more comprehensive and realistic 
representation of the material response under varying thermal 
and loading conditions. One of the drawbacks we faced was the 
time consumption and the difficulty of providing practical 
examples without resorting to numerical analysis and simulation.}

{We have numerous prospects ahead. In particular, 
we intend to conduct work on numerical analysis 
and simulation based on the results obtained 
for the studied problem.}


\vspace{6pt} 

\authorcontributions{Conceptualization, B.F. and M.D.; 
validation, B.F., M.D. and D.F.M.T.; 
formal analysis, B.F., M.D. and D.F.M.T.; 
investigation, B.F., M.D. and D.F.M.T.; 
writing---original draft preparation, B.F., M.D. and D.F.M.T.; 
writing---review and editing, B.F., M.D. and D.F.M.T.; 
supervision, M.D. and D.F.M.T.; 
project administration, M.D. and D.F.M.T.; 
funding acquisition, B.F., M.D. and D.F.M.T. 
All authors have read and agreed to the published version of the manuscript.}

\funding{This research was funded by 
Fundação para a Ciência e a Tecnologia 
grant number UIDB/04106/2020 
(\url{https://doi.org/10.54499/UIDB/04106/2020}, 
accessed on 25 January 2024).} 

\dataavailability{Data are contained within the article.} 

\acknowledgments{Fadlia is grateful for the financial support 
of University of Constantine 1, Algeria, 
for a one-month visit to the R\&D Unit CIDMA, 
Department of Mathematics, University of Aveiro. 
The hospitality of the host institution 
is here gratefully acknowledged. 
{The authors are grateful to four anonymous
referees for several constructive comments and remarks.}}

\conflictsofinterest{The authors declare no conflicts of interest.} 


\begin{adjustwidth}{-\extralength}{0cm}
	
\reftitle{References}


\PublishersNote{}

\end{adjustwidth}

\end{document}